\documentclass[11pt]{article}
\usepackage{amssymb}

\begin{document}



\newcommand{\newc}{\newcommand}


\renewcommand{\theequation}{\thesection.\arabic{equation}}
\newc{\eqnoset}{\setcounter{equation}{0}}
\newcommand{\myref}[2]{#1~\ref{#2}}

\newcommand{\mref}[1]{(\ref{#1})}
\newcommand{\reflemm}[1]{Lemma~\ref{#1}}
\newcommand{\refrem}[1]{Remark~\ref{#1}}
\newcommand{\reftheo}[1]{Theorem~\ref{#1}}
\newcommand{\refdef}[1]{Definition~\ref{#1}}
\newcommand{\refcoro}[1]{Corollary~\ref{#1}}
\newcommand{\refprop}[1]{Proposition~\ref{#1}}
\newcommand{\refsec}[1]{Section~\ref{#1}}
\newcommand{\refchap}[1]{Chapter~\ref{#1}}

\newcommand{\beq}{\begin{equation}}
\newcommand{\eeq}{\end{equation}}
\newcommand{\beqno}[1]{\begin{equation}\label{#1}}

\newcommand{\barr}{\begin{array}}
\newcommand{\earr}{\end{array}}

\newc{\bearr}{\begin{eqnarray*}}
\newc{\eearr}{\end{eqnarray*}}

\newc{\bearrno}[1]{\begin{eqnarray}\label{#1}}
\newc{\eearrno}{\end{eqnarray}}

\newc{\non}{\nonumber}
\newc{\nol}{\nonumber\nl}

\newcommand{\bdes}{\begin{description}}
\newcommand{\edes}{\end{description}}
\newc{\benu}{\begin{enumerate}}
\newc{\eenu}{\end{enumerate}}
\newc{\btab}{\begin{tabular}}
\newc{\etab}{\end{tabular}}



\newtheorem{theorem}{Theorem}[section]
\newtheorem{defi}[theorem]{Definition}
\newtheorem{lemma}[theorem]{Lemma}
\newtheorem{rem}[theorem]{Remark}
\newtheorem{exam}[theorem]{Example}
\newtheorem{propo}[theorem]{Proposition}
\newtheorem{corol}[theorem]{Corollary}

\renewcommand{\thelemma}{\thesection.\arabic{lemma}}

\newcommand{\btheo}[1]{\begin{theorem}\label{#1}}
\newc{\brem}[1]{\begin{rem}\label{#1}\em}
\newc{\bexam}[1]{\begin{exam}\label{#1}\em}
\newc{\bdefi}[1]{\begin{defi}\label{#1}}
\newcommand{\blemm}[1]{\begin{lemma}\label{#1}}
\newcommand{\bprop}[1]{\begin{propo}\label{#1}}
\newcommand{\bcoro}[1]{\begin{corol}\label{#1}}
\newcommand{\etheo}{\end{theorem}}
\newcommand{\elemm}{\end{lemma}}
\newcommand{\eprop}{\end{propo}}
\newcommand{\ecoro}{\end{corol}}
\newc{\erem}{\end{rem}}
\newc{\eexam}{\end{exam}}
\newc{\edefi}{\end{defi}}

\newc{\rmk}[1]{{\bf REMARK #1: }}
\newc{\DN}[1]{{\bf DEFINITION #1: }}

\newcommand{\bproof}{{\bf Proof:~~}}
\newc{\eproof}{{\vrule height8pt width5pt depth0pt}\vspace{3mm}}


\newcommand{\rarr}{\rightarrow}
\newcommand{\Rarr}{\Rightarrow}
\newcommand{\tru}{\backslash}
\newc{\bfrac}[2]{\dspl{\frac{#1}{#2}}}


\newc{\nl}{\vspace{2mm}\\}
\newc{\nid}{\noindent}


\newcommand{\oneon}[1]{\frac{1}{#1}}
\newcommand{\dspl}{\displaystyle}
\newc{\grad}{\nabla}
\newc{\Div}{\mbox{div}}
\newc{\pdt}[1]{\dspl{\frac{\partial{#1}}{\partial t}}}
\newc{\pdn}[1]{\dspl{\frac{\partial{#1}}{\partial \nu}}}
\newc{\pdNi}[1]{\dspl{\frac{\partial{#1}}{\partial \mathcal{N}_i}}}
\newc{\pD}[2]{\dspl{\frac{\partial{#1}}{\partial #2}}}
\newc{\dt}{\dspl{\frac{d}{dt}}}
\newc{\bdry}[1]{\mbox{$\partial #1$}}
\newc{\sgn}{\mbox{sign}}

\newc{\Hess}[1]{\frac{\partial^2 #1}{\pdh z_i \pdh z_j}}
\newc{\hess}[1]{\partial^2 #1/\pdh z_i \pdh z_j}


\newcommand{\Coone}[1]{\mbox{$C^{1}_{0}(#1)$}}
\newcommand{\lspac}[2]{\mbox{$L^{#1}(#2)$}}
\newc{\hspac}[2]{\mbox{$C^{0,#1}(#2)$}}
\newc{\Hspac}[2]{\mbox{$C^{1,#1}(#2)$}}
\newc{\Hosp}{\mbox{$H^{1}_{0}$}}
\newcommand{\Lsp}[1]{\mbox{$L^{#1}(\Og)$}}
\newc{\hsp}{\Hosp(\Og)}


\newc{\ag}{\alpha}
\newc{\bg}{\beta}
\newc{\cg}{\gamma}\newc{\Cg}{\Gamma}
\newc{\dg}{\delta}\newc{\Dg}{\Delta}
\newc{\eg}{\varepsilon}
\newc{\zg}{\zeta}
\newc{\thg}{\theta}
\newc{\llg}{\lambda}\newc{\LLg}{\Lambda}
\newc{\kg}{\kappa}
\newc{\rg}{\rho}
\newc{\sg}{\sigma}\newc{\Sg}{\Sigma}
\newc{\tg}{\tau}
\newc{\fg}{\phi}\newc{\Fg}{\Phi}
\newc{\vfg}{\varphi}
\newc{\og}{\omega}\newc{\Og}{\Omega}
\newc{\pdh}{\partial}

\newc{\ccG}{{\cal G}}


\newc{\ii}[1]{\int_{#1}}
\newc{\iidx}[2]{{\dspl\int_{#1}~#2~dx}}
\newc{\bii}[1]{{\dspl \ii{#1} }}
\newc{\biii}[2]{{\dspl \iii{#1}{#2} }}
\newc{\su}[2]{\sum_{#1}^{#2}}
\newc{\bsu}[2]{{\dspl \su{#1}{#2} }}

\newcommand{\iiomdx}[1]{{\dspl\int_{\Og}~ #1 ~dx}}
\newc{\biiom}[1]{{\dspl\int_{\bdrom}~ #1 ~d\sg}}
\newc{\io}[1]{{\dspl\int_{\Og}~ #1 ~dx}}
\newc{\bio}[1]{{\dspl\int_{\bdrom}~ #1 ~d\sg}}
\newc{\bsir}{\bsu{i=1}{r}}
\newc{\bsim}{\bsu{i=1}{m}}

\newc{\iibr}[2]{\iidx{\bprw{#1}}{#2}}
\newc{\Intbr}[1]{\iibr{R}{#1}}
\newc{\intbr}[1]{\iibr{\rg}{#1}}
\newc{\intt}[3]{\int_{#1}^{#2}\int_\Og~#3~dxdt}

\newc{\itQ}[2]{\dspl{\int\hspace{-2.5mm}\int_{#1}~#2~dz}}
\newc{\mitQ}[2]{\dspl{\rule[1mm]{4mm}{.3mm}\hspace{-5.3mm}\int\hspace{-2.5mm}\int_{#1}~#2~dz}}
\newc{\mitQQ}[3]{\dspl{\rule[1mm]{4mm}{.3mm}\hspace{-5.3mm}\int\hspace{-2.5mm}\int_{#1}~#2~#3}}

\newc{\mitx}[2]{\dspl{\rule[1mm]{3mm}{.3mm}\hspace{-4mm}\int_{#1}~#2~dx}}
\newc{\mitmu}[2]{\dspl{\rule[1mm]{3mm}{.3mm}\hspace{-4mm}\int_{#1}~#2~d\mu}}
\newc{\iidmu}[2]{{\dspl\int_{#1}~#2~d\mu}}

\newc{\mitQq}[2]{\dspl{\rule[1mm]{4mm}{.3mm}\hspace{-5.3mm}\int\hspace{-2.5mm}\int_{#1}~#2~d\bar{z}}}
\newc{\itQq}[2]{\dspl{\int\hspace{-2.5mm}\int_{#1}~#2~d\bar{z}}}

\newc{\pder}[2]{\dspl{\frac{\partial #1}{\partial #2}}}


\newc{\ui}{u_{i}}
\newcommand{\upl}{u^{+}}
\newcommand{\umn}{u^{-}}
\newcommand{\un}{\{ u_{n}\}}

\newcommand{\uo}{u_{0}}
\newc{\voi}{v_{i}^{0}}
\newc{\uoi}{u_{i}^{0}}
\newc{\vu}{\vec{u}}

\newc{\xo}{x_{0}}
\newc{\Br}{B_{R}}
\newc{\Bro}{\Br (\xo)}
\newc{\bdrom}{\bdry{\Og}}
\newc{\ogr}[1]{\Og_{#1}}
\newc{\Bxo}{B_{x_0}}

\newc{\inP}[2]{\|#1(\bullet,t)\|_#2\in\cP}
\newc{\cO}{{\mathcal O}}
\newc{\inO}[2]{\|#1(\bullet,t)\|_#2\in\cO}

\newc{\newl}{\\ &&}

\newc{\bilhom}{\mbox{Bil}(\mbox{Hom}(\RR^{nm},\RR^{nm}))}
\newc{\VV}[1]{{V(Q_{#1})}}

\newc{\ccA}{{\mathcal A}}
\newc{\ccB}{{\mathcal B}}
\newc{\ccC}{{\mathcal C}}
\newc{\ccD}{{\mathcal D}}
\newc{\ccE}{{\mathcal E}}
\newc{\ccH}{\mathcal{H}}
\newc{\ccF}{\mathcal{F}}
\newc{\ccI}{{\mathcal I}}
\newc{\ccJ}{{\mathcal J}}
\newc{\ccP}{{\mathcal P}}
\newc{\ccQ}{{\mathcal Q}}
\newc{\ccR}{{\mathcal R}}
\newc{\ccS}{{\mathcal S}}
\newc{\ccT}{{\mathcal T}}
\newc{\ccX}{{\mathcal X}}
\newc{\ccY}{{\mathcal Y}}
\newc{\ccZ}{{\mathcal Z}}

\newc{\bb}[1]{{\mathbf #1}}
\newc{\bbA}{{\mathbf A}}
\newc{\myprod}[1]{\langle #1 \rangle}
\newc{\mypar}[1]{\left( #1 \right)}

\newc{\BLLg}{\mathbf{\LLg}}

\newc{\mA}{\mathbf{A}}
\newc{\mB}{\mathbf{B}}
\newc{\mC}{\mathbf{C}}
\newc{\mD}{\mathbf{D}}
\newc{\mE}{\mathbf{E}}
\newc{\mF}{\mathbf{F}}
\newc{\mJ}{\mathbf{J}}
\newc{\mG}{\mathbf{G}}
\newc{\mP}{\mathbf{P}}
\newc{\mR}{\mathbf{R}}
\newc{\mQ}{\mathbf{Q}}
\newc{\mX}{\mathbf{X}}
\newc{\muu}{\mathbf{u}}
\newc{\mvv}{\mathbf{v}}

\newc{\mllg}{\mathbb{\lambda}}
\newc{\mLLg}{\mathbf{\LLg}}


\newc{\lspn}[2]{\mbox{$\| #1\|_{\Lsp{#2}}$}}
\newc{\Lpn}[2]{\mbox{$\| #1\|_{#2}$}}
\newc{\Hn}[1]{\mbox{$\| #1\|_{H^1(\Og)}$}}


\newc{\cyl}[1]{\og\times \{#1\}}
\newc{\cyll}{\og\times[0,1]}
\newc{\vx}[1]{v\cdot #1}
\newc{\vtx}[1]{v(t,x)\cdot #1}
\newc{\vn}{\vx{n}}

\newcommand{\RR}{{\rm I\kern -1.6pt{\rm R}}}


\newenvironment{proof}{\noindent\textbf{Proof.}\ }
{\nopagebreak\hbox{ }\hfill$\Box$\bigskip}


\newc{\itQQ}[2]{\dspl{\int_{#1}#2\,dz}}
\newc{\mmitQQ}[2]{\dspl{\rule[1mm]{4mm}{.3mm}\hspace{-4.3mm}\int_{#1}~#2~dz}}
\newc{\MmitQQ}[2]{\dspl{\rule[1mm]{4mm}{.3mm}\hspace{-4.3mm}\int_{#1}~#2~d\mu}}

\newc{\MUmitQQ}[3]{\dspl{\rule[1mm]{4mm}{.3mm}\hspace{-4.3mm}\int_{#1}~#2~d#3}}
\newc{\MUitQQ}[3]{\dspl{\int_{#1}~#2~d#3}}

\newc{\mccP}{\mathbb{P}}
\newc{\mccK}{\mathbb{K}}

\newc{\DKTmU}{\mccK(U)}
\newc{\DKTmUold}{(K_U(U)^{-1})^T}

\newc{\myPi}{\mathbf{W}}
\newc{\myIbar}{\bar{\ccI}_1}
\newc{\myIhat}{\hat{\ccI}_1}
\newc{\myIbreve}{\breve{\ccI}_0}

\vspace*{-.8in}
\begin{center} {\LARGE\em Weighted Gagliardo-Nirenberg Inequalities Involving BMO Norms and Measures}

 \end{center}

\vspace{.1in}

\begin{center}

{\sc Dung Le}{\footnote {Department of Mathematics, University of
Texas at San
Antonio, One UTSA Circle, San Antonio, TX 78249. {\tt Email: Dung.Le@utsa.edu}
}}

\end{center}

\begin{abstract}
Global and local weighted Gagliardo-Nirenberg inequalities with doubling measures are established. These inequalities are key ingredients for the regularity theory and existence of strong solutions for strongly coupled parabolic and elliptic systems which are degenerate or singular because of the unboundedness of dependent and independent variables. \end{abstract}

\vspace{.2in}

\section{Introduction}\label{introsec}\eqnoset

In \cite{SRR,SR},  for any $p\ge1$ and $C^2$ scalar function $u$ on $\RR^n$, $n\ge2$, global and local Gagliardo-Nirenberg inequalities of the form \beqno{GNSR}\iidx{\RR^n}{|Du|^{2p+2}} \le C(n,p)\|u\|_{BMO}^2\iidx{\RR^n}{|Du|^{2p-2}|D^2u|^2}\eeq were established and applied to the solvability of {\em scalar} elliptic equations.

More general and vectorial versions of these inequalities were presented in \cite{letrans,dleANS} to establish the solvability of strongly counpled parabolic systems of the form {\em  nonregular} but {\em uniform}  parabolic system
\beqno{e10}\left\{\barr{ll} u_t=\Div(A(u,Du))+\hat{f}(u,Du)& (x,t)\in Q=\Og\times(0,T_0),\\u(x,0)=U_0(x)& x\in\Og\\\mbox{$u=0$ on $\partial \Og\times(0,T_0)$}. &\earr\right.\eeq

Here, and throughout this paper, $\Og$ is a bounded domain with smooth boundary $\partial \Og$ in $\RR^d$ for some integer $d\ge1$. The temporal and $k$-order spatial derivatives of a vector-valued function $$u(x,t)=(u_1(x,t),\ldots,u_m(x,t))^T \quad m>1 $$ are denoted by $u_t$ and $D^ku$ respectively.

In this paper, we generalize global and local versions of \mref{GNSR} and the inequalities in \cite{dleANS} (see \refcoro{GNlobalANS1mcoro}).  
Roughly speaking, we will establish inequalities of the following type: for any $p\ge1$ and any  $C^2$ map $U:\Og\to\RR^m$
$$\barr{ll}\lefteqn{\iidmu{\Og}{\Fg^2(U)|DU|^{2p+2}}\le}\hspace{2cm}&\\& C\|K(U)\|_{BMO(\mu)}^2\{\iidmu{\Og}{\LLg^2(U)|DU|^{2p-2}|D^2U|^2}+\cdots\}.\earr$$
Here, $K$ is a map, $\Fg,\LLg$ are functions on $\RR^m$ and $d\mu=\og dx$ is a doubling measure on $\Og$. We assume that $\Og, \mu$ support a Poincar\'e-Sobolev type inequality.

The purpose of such generalization becomes clear when we apply the results to the study of  local/global existence of strong solutions to \mref{e10} in our forthcoming work \cite{dlepde}.
First of all, by replacing the Lebesgue measure $dx$ with a general measure $d\mu=\og dx$, we allow the matrices $A,\hat{f}$ in \mref{e10} to depend on $x,t$ and become degenerate or singular near a subset of $\Og$.

Secondly, the degeneracy and singularity of \mref{e10} can also come from the behavior of the solution $u$ itself, which is not well known as maximum principles are not available for systems (i.e. $m>1$). We replace the factor $\|u\|_{BMO}$ in \mref{GNSR} by $\|K(u)\|_{BMO(\mu)}$ where $K$ is a map in $\RR^m$. This allows us to deal with the case when estimates for $\|u\|_{BMO}$, but $\|K(u)\|_{BMO(\mu)}$, are not available. For example, one of the consequences of our general inequalities in this paper is the following inequality which will be useful in dealing with degenerate system in \cite{dlepde}: if $\|\log(|u|)\|_{BMO(\mu)}$ is sufficiently small then 
$$\barr{ll}\lefteqn{\iidmu{\Og}{|u|^{2k-2}|Du|^{2p+2}}\le}\hspace{1cm}&\\& C\|\log(|u|)\|_{BMO(\mu)}^2\iidmu{\Og}{(|u|^{2k}|Du|^{2p-2}|D^2u|^2+|u|^{2k}|Du|^{2p})}.\earr$$

Various choices of $K$ will be discussed in \cite{dlepde}.

We organize the paper as follows. The hypotheses and main results will be presented in \refsec{mainres}. One of our key ingredients of the proof comes from Tolsa's work \cite{tolhh} on the $RBMO(\mu)$ spaces and we will discuss it in \refsec{tolsasec}. The main global and technical inequality is stated and proved in \refsec{GNineqsec}. The local version is then established in \refsec{localsec}. We conclude the paper with the proof of the main inequalities and their consequences in \refsec{collsec}.

\section{Hypotheses and Main results}\eqnoset\label{mainres}

Throughout this paper, in  our statements and proofs, we use $C,C_1,\ldots$ to denote various constants which can change from line to line but depend only on the parameters of the hypotheses in an obvious way. We will write $C(a,b,\ldots)$ when the dependence of  a constant $C$ on its parameters is needed to emphasize that $C$ is bounded in terms of its parameters. We also write $a\preceq b$ if there is a universal constant $C$ such that $a\le Cb$. In the same way, $a\sim b$ means $a\preceq b$ and $b\preceq a$.

For any $\mu$-measurable subset $A$ of $\Og$  and any  locally $\mu$-integrable function $U:\Og\to\RR^m$ we denote by  $\mu(A)$ the measure of $A$ and $U_A$ the average of $U$ over $A$. That is, $$U_A=\mitmu{A}{U(x)} =\frac{1}{\mu(A)}\iidmu{A}{U(x)}.$$

We say that $\Og$ and $\mu$ support a Poincar\'e-Sobolev inequality if the following holds. \bdes \item[PS)]  There are $\sg\in(0,1)$ and $\tau_*\ge 1$ such that for {\em some} $q>2$ and $q_*=\sg q<2$ we have
\beqno{PSineq} \frac{1}{l(B)}
\left(\mitmu{B}{|u-u_{B}|^q}\right)^\frac1q \le C_{PS}
\left(\mitmu{\tau_*B}{|Du|^{q_*}}\right)^\frac{1}{q_*}\eeq for some constant $C_{PS}$ and any cube $B$ with side length $l(B)$ and any function $u\in C^1(B)$. \edes

Here and throughout this paper, we write $B_R(x)$ for a cube centered at $x$ with side length $R$ and sides parallel to to standard axes of $\RR^d$. We will omit $x$ in the notation $B_R(x)$ if no ambiguity can arise. We denote by $l(B)$ the side length of $B$ and by $\tau B$ the cube which is concentric with $B$ and has side length $\tau l(B)$.

We have the following remark on the validity of the assumption PS).
\brem{PSrem}
Suppose that  $\mu$ is doubling and supports a $q_*$-Poincar\'e inequality (see \cite[eqn. (5)]{Haj}): There are some constants $C_P$, $q_*\in[1,2]$ and $\tau_*\ge1$ the following inequality holds true
\beqno{Pineq} 
\mitmu{B}{|h-h_{B}|}\le C_Pl(B)
\left(\mitmu{\tau_*B}{|Dh|^{q_*}}\right)^\frac{1}{q_*}\eeq for any cube $B$ with side length $l(B)$ and any function $h\in C^1(B)$.

Assume also that for some $s>0$ $\mu$ satisfies  the following inequality \beqno{fracmu} \left(\frac{r}{r_0}\right)^s\preceq\frac{\mu(B_r(x))}{\mu(B_{r_0}(x_0))},\eeq where $B_r(x), B_{r_0}(x_0)$ are any cubes with $x\in B_{r_0}(x_0)$. If $q_*=2$ then \cite[Section 3]{Haj} shows that  a $q_*$-Poincar\'e inequality also holds for some $q_*<2$. Thus, we can assume that $q_*\in(1,2)$.

If $q_*<s$ then \cite[1) of Theorem 5.1]{Haj} establishes \mref{PSineq}
for  $q=sq_*/(s-q_*)$. Thus, $q>2$ if $s<2q_*/(2-q_*)$. This is the case if we choose $q_*<2$ and closed to 2.  If $s=q_*$, \cite[2) of Theorem 5.1]{Haj} shows that \mref{PSineq} holds true
for any $q>1$. On the other hand, if $q_*>s$ then \cite[3) of Theorem 5.1]{Haj} gives a stronger version of \mref{PSineq} for $q=\infty$. In particular, the H\"older norm of $u$ is bounded in terms of $\|Du\|_{L^{q_*}(\mu)}$. We thus need only that  $\Og,\mu$ support a $q_*$-Poincar\'e inequality for some $q_*\in(1,2)$ and \mref{fracmu} is valid for some $s>0$.

\erem

To proceed, we recall some well known notions from Harmonic Analysis.

A function $f\in L^1(\mu)$ is said to be in $BMO(\mu)$ if \beqno{bmodef} [f]_{*,\mu}:=\sup_Q\mitmu{Q}{|f-f_Q|}<\infty.\eeq We then define $$\|f\|_{BMO(\mu)}:=[f]_{*,\mu}+\|f\|_{L^1(\mu)}.$$

For $\cg\in(1,\infty)$ we say that a nonnegative locally integrable function $w$ belongs to the class $A_\cg$ or $w$ is an $A_\cg$ weight if the quantity
\beqno{aweight} [w]_{\cg} := \sup_{B\subset\Og} \left(\mitmu{B}{w}\right) \left(\mitmu{B}{w^{1-\cg'}}\right)^{\cg-1} \quad\mbox{is finite}.\eeq
Here, $\cg'=\cg/(\cg-1)$ and the supremum is taken over all cubes $B$ in $\Og$. For more details on these classes we refer the reader to \cite{FPW,OP,st}.

We assume the following hypotheses.
\bdes \item[M)] Let $\Og$ be a bounded domain in $\RR^d$ and $d\mu=\og dx$ for some $\og\in C^1(\Og,\RR^+)$. Suppose that there are  a constant $C_\mu$ and a fixed number $n\in(0,d]$ such that : for any cube $Q_r$ with side length $r>0$ \beqno{sizecond1m}\mu(Q_r)\le C_\mu r^n.\eeq Furthermore, $\Og,\mu$ satisfies PS).

\item[A.1)] Let $K:\mbox{dom}(K)\to\RR^m$ be a $C^1$ map on a domain $\mbox{dom}(K)\subset\RR^m$ such that $K_U^{-1}(U)=K_U(U)^{-1}$ exists and $\mccK_U\in L^\infty(\mbox{dom}K)$, where we will always abbreviate \beqno{DKdef} \DKTmU=(K_U(U)^{-1})^T.\eeq

Let $\Fg,\LLg:\mbox{dom}(K)\to\RR^+$ be $C^1$ positive functions. We assume that for all $U\in \mbox{dom}(K)$
\beqno{kappamain} |\DKTmU|\preceq \LLg(U)\Fg^{-1}(U),\eeq
\beqno{logcondszmain}|\Fg_U(U)||\mccK(U)|\preceq \Fg(U) \mbox{ and } |\mccK_U(U)| \mbox{ is bounded}.\eeq
\edes

\bdes

\item[A.2)] Let $U:\Og\to \mbox{dom}(K)$ be a $C^2$ map such that
\beqno{boundarym}\myprod{\og\Fg^2(U)\DKTmU DU,\vec{\nu}}=0\eeq on $\partial\Og$ where $\vec{\nu}$ is the outward normal vector of $\partial\Og$.

\item[A.3)] Let $\myPi(x):=\LLg^{p+1}(U(x))\Fg^{-p}(U(x))$. Asume that $[\myPi^{\ag}]_{\bg+1}$ is finite for some $\ag>2/(p+2)$ and $\bg<p/(p+2)$.
\edes

We denote
\beqno{Idefm} I_1:=\iidmu{\Og}{\Fg^2(U)|DU|^{2p+2}},\;
I_2:=\iidmu{\Og}{\LLg^2(U)|DU|^{2p-2}|D^2U|^2},\eeq
\beqno{Idef1zm} \myIbar:=\iidmu{\Og}{|\LLg_U(U)|^2|DU|^{2p+2}},\eeq
\beqno{Idef1zmm}\myIbreve:=\iidmu{\Og}{|D\og|^2\og^{-2}\LLg^2(U)|DU|^{2p}}.\eeq

Our first main result is the following.

\btheo{GNlobalANSm} Assume A.1)-A.3). There are constants $C,C([\myPi^{\ag}]_{\bg+1})$  for which \beqno{GNglobalestzm}I_1\le C\|K(U)\|_{BMO}^2\left[I_2+\myIbar+C([\myPi^{\ag}]_{\bg+1}) [I_2+I_1+\myIbreve]\right].\eeq 

In addition, if \beqno{KeyLLgPicond1m}|\LLg_U|\preceq \Fg.\eeq
Then there is a constant $C([\myPi^{\ag}]_{\bg+1})$  such that \beqno{GNglobalestam}I_1\le C\|K(U)\|_{BMO}^2\left[I_2+I_1+C([\myPi^{\ag}]_{\bg+1}) [I_2+I_1+\myIbreve]\right].\eeq
Here, $C$ also depends on $C_{PS},C_\mu$.

\etheo

Next, we have a local version of \mref{GNglobalest}. Let $\Og_*$ be a subset of $\Og$. In place of M), the condition on the measure $d\mu$, we assume that there are two functions $\og_*,\og_0$ satisfying the following conditions.

\bdes \item[LM.0)] 
$\og_*\in C^1(\Og)$ and satisfies $\og_*\equiv 1$ in $\Og_*$ and $\og_*\le 1$ in $\Og$.
\beqno{subogm} \og_*\equiv 1 \mbox{ in $\Og_*$ and } \og_*\le 1 \mbox{ in $\Og$}.\eeq

\item[LM.1)] $\og_0\in C^1(\Og)$ and for $d\mu=\og_0^2dx$ and some $n\in(0,d]$ we have $\mu(B_r)\le Cr^n$. 

\item[LM.2)]  $d\mu=\og_0^2 dx$ supports the Poincar\'e-Sobolev inequality \mref{PSineq} in PS). In addition, $\og_0$ also supports a Hardy type inequality: There is a constant $C_H$ such that for any function $u\in C^1_0(B)$\beqno{lehr1m} \iidx{\Og}{|u|^2|D\og_0|^2}\le C_H\iidx{\Og}{|Du|^2\og_0^2}\eeq 
\edes

\btheo{GNlocalog1m} Suppose LM.0)-LM.2), A.1)-A.3) and that (compare to \mref{boundarym} with $\og$ being $\og_*\og_0^2$)
\beqno{boundaryzm}\myprod{\og_*\og_0^2 \Fg^2(U)\DKTmU DU,\vec{\nu}}=0\eeq on $\partial\Og$ where $\vec{\nu}$ is the outward normal vector of $\partial\Og$.

For any $\og_1\in L^1(\Og)$ and $\og_1\sim \og_0^2$ we define $d\mu=\og_1 dx$ and \beqno{I1*m}I_{1,*}:=\iidmu{\Og_*}{\Fg^2(U)|DU|^{2p+2}},\eeq \beqno{I0*m} \breve{I}_{0,*}:=\sup_\Og|D\og_*|^2\iidmu{\Og}{\LLg^2(U)|DU|^{2p}}.\eeq

Then, for any $\eg>0$ there are constants $C, C([\myPi^{\ag}]_{\bg+1})$ such that
\beqno{GNlocog11m}I_{1,*}\le  \eg I_1+\eg^{-1}C\|K(U)\|_{BMO(\mu)}^2[I_2+\myIbar+C([\myPi^{\ag}]_{\bg+1}) [I_2+\myIbar+\breve{I}_{0,*}]].\eeq
Here, $C$ also depends on $C_{PS},C_\mu$ and $C_H$.

\etheo

\brem{dogrem} A typical choice of $\og_0$ that satisfies the Hardy type inequality \mref{lehr1m} in LM.2)  is $\og_0(x)=d_\Og^\frac{\cg}{2}(x)$. Then, for $d\mu=\og_1 dx\sim d_\Og^\cg dx$ we will check the conditions LM.1)-LM.2). If $B_r$ is far away from $\partial\Og$, we have $\mu(B_r)\preceq r^d$. Near the boundary, as $\partial\Og$ is $C^1$, we easily see that $\mu(B_r)\preceq   r^{d+\bg}$. If $d+\cg\ge n$ and $r$ is bounded then $\mu(B_r)\le Cr^n$. This is the case because $\Og$ is bounded.  Thus, for any $\cg> -d$, we define $n=\min\{d,d+\cg\}\in(0,d]$ to see that $\mu(B_r)\le Cr^n$ for some constant $C$ which is bounded in terms of $\mbox{diam}(\Og)$.

We now recall the following Hardy inequality proved by Necas (see also the paper by Lehrb\"ack \cite{Lehr} for much more general versions)
\beqno{lehr} \iidx{\Og}{|u(x)|^qd_\Og^{\cg-q}(x)}\le C\iidx{\Og}{|Du(x)|^qd_\Og^{\cg}(x)},\quad \cg<q-1.\eeq

We see that \mref{lehr} in LM.2) holds true with $q=2$, $\cg<1$.
\erem

An immediate consequence of \reftheo{GNlobalANSm} is the following main inequality in \cite{dleANS}.

\bcoro{GNlobalANS1mcoro} Let $U:\Og\to \mbox{dom}(K)$ be a $C^2$ vector-valued function. Suppose that either $U$ or $\Fg^2(u)\frac{\partial U}{\partial \nu}$ vanish on the boundary $\partial \Og$  of $\Og$. 

We set \beqno{Idefaz} I_1:=\iidx{\Og}{\Fg^2(U)|DU|^{2p+2}},\; I_2:=\iidx{\Og}{\Fg^2(U)|DU|^{2p-2}|D^2U|^2},\eeq
\beqno{Idefazz} \myIbar:=\iidx{\Og}{|\Fg_U(U)|^2|DU|^{2p+2}}.\eeq

For any $\ag>2/(p+2)$ and $\bg<p/(p+2)$ we have \beqno{GNglobalestazz}I_1 \le  C\|U\|_{BMO(\Og)}^2
\left[I_2+\myIbar + C([\Fg^{\ag}(U)]_{\bg+1})(I_2+I_1+\myIbar)\right].\eeq
\ecoro

{\bf Proof of \refcoro{GNlobalANS1mcoro}:}
We simply choose $\Fg=\LLg$ and $K(U)=U$ to see $\myPi:=\Fg$. For $\og\equiv1$, $\mu$ is then the Lebesgue measure. As $\myIbar$ in \mref{Idef1zm} and \mref{Idefazz} are the same  and $\myIbreve=0$, we then have from \mref{GNglobalestazz} from \mref{GNglobalestzm}. \eproof

\reftheo{GNlocalog1m} with $\og_0\equiv1$ also implies the local version of \refcoro{GNlobalANS1mcoro} which is one of the key ingredients in the proof of solvability of strongly coupled parabolic systems in \cite{dleANS}. In this paper, we obtain a more general result with general $\mu$ satisfying LM.0)-LM.2). \beqno{GNlocog11ANS}I_{1,*}\le  \eg I_1+\eg^{-1}C\|U\|_{BMO(\Og)}^2[I_2+\myIbar+C([\Fg^{\ag}]_{\bg+1}) [I_2+\myIbar+\breve{I}_{0,*}]].\eeq

Of course, there are many ways to choose $K,\LLg,\Fg$ depending on different situations in applications. Let us consider another choice of $K$ and the connection between the two terms $\|K(U)\|_{BMO(\mu)}$, $[\myPi^{\ag}]_{\bg+1}$. In this paper we will only look at the case  $K(U)=[\log(\eg+|U_i|)]_{i=1}^m$, which will be useful in dealing with porous media type parabolic systems in our forthcoming work \cite{dlepde}. Different choices of $K$ will be presented in \cite{dlepde} too.

For $\LLg(U)=(\eg+|U|)^k$ and $\Fg(U)\sim|\LLg_U(U)|$ we then define for any $k\ne0$ and $\eg\ge0$
\beqno{logI1m}I_1= \iidmu{\Og}{(\eg+|U|)^{2k-2}|DU|^{2p+2}},\;\myIbreve= \iidmu{\Og}{(\eg+|U|)^{2k}|DU|^{2p}},\eeq
\beqno{logI2m}I_2= \iidmu{\Og}{(\eg+|U|)^{2k}|DU|^{2p-2}|D^2U|^2}.\eeq

\bcoro{GNlobalANS} For $m\ge1$, any $k\ne0$ and $\eg\ge0$ we consider the map \beqno{Klogm}K(U)=[\log(\eg+|U_i|)]_{i=1}^m, \quad U=[U_i]_{i=1}^m.\eeq With the notations \mref{logI1m} and \mref{logI2m} and $\myPi=(\eg+|U|)^{k+p}$, we have
\beqno{GNglobalestalogm}I_1\le C\|K(U)\|_{BMO(\mu)}^2\left[I_2+I_1+C([\myPi^{\ag}]_{\bg+1}) [I_2+I_1+\myIbreve]\right],\eeq as long as the integrals are finite. Here, $C$ is independent of $\eg$.
\ecoro

\newc{\mmc}{\mathbf{c}}

We consider the case $m=1$. As $\myPi=\LLg^{p+1}\Fg^{-p}=|k|^{-p}(\eg+|U|)^{k+p}$,
we have $[\myPi^\ag]_{A_q}=[(\eg+|U|)^{\ag(k+p)}]_{A_q}$  and $$[\log(\myPi^\ag)]_{*,\mu}=\ag|k+p|[\log(\eg+|U|)]_{*,\mu}.$$ 

Via  a simple use of Jensen's inequality, it is well known (e.g. see \cite[Chapter 9]{Graf}) that $\|\log w\|_{BMO}\le [w]_{A_q}$ for $1<q\le2$. In our case, $q=\bg+1<2$ so that  $\|\log\myPi^\ag\|_{BMO}\le [\myPi^\ag]_{A_q}$. Thus, the term $[\log(\eg+|U|)]_{BMO(\mu)}$ can be controlled by $[\myPi^\ag]_{A_q}$. However, this type of result is not helpful in the regularity theory of PDEs.

On the other hand, if $\log \myPi$ is BMO then we also know that $\myPi$ is a weight. We recall the following John-Nirenberg inequality (e.g. see \cite[Chapter 9]{Graf}): If $\mu$ is doubling then for any BMO($\mu$) function $v$ 
there are constants $\mmc_1,\mmc_2$, which depend only on the doubling constant of $\mu$, such that  
\beqno{JNineqm}\mitmu{B}{e^{\frac{\mmc_1}{[v]_{*,\mu}}|v-v_B|}} \le \mmc_2.\eeq

We then have the following result.

\bcoro{GNlobalANS1} In addition to the assumptions of \refcoro{GNlobalANS} we suppose that \beqno{Westcond1} |k+p|[\log(\eg +|U|)]_{*,\mu}\le \mmc_1\bg\ag^{-1}.\eeq

Then there is a constant $C$, which depends also on $\mmc_2$, for which \beqno{GNglobalestzzz}I_1\le C\|\log(\eg+|U|)\|_{BMO(\mu)}^2\left[I_2+I_1+\myIbreve\right].\eeq
\ecoro

It is clear that if $\|\log(\eg+|U|)\|_{BMO(\mu)}$ is sufficiently small then \mref{Westcond1} and \mref{GNglobalestzzz} imply $$I_1\le C\|\log(\eg+|U|)\|_{BMO(\mu)}^2\left[I_2+\myIbreve\right].$$

Of course, the above corollaries have their local versions from \reftheo{GNlocalog1m}.

\section{Some simple consequences from Tolsa's works}\eqnoset\label{tolsasec}

The $RBMO(\mu)$ space was introduced by Tolsa in \cite{tol3,tolhh}. Tolsa considered {\em non-doubling} measure $\mu$ and defined \beqno{rbmodef}[f]_{*,\mu}:=\sup_Q\mitmu{\llg Q}{|f-f_Q|}\eeq for some constant $\llg>1$. This constant $\llg$ is not important as shown in \cite{tolhh}. The definition of $RBMO(\mu)$ spaces in \cite{tolhh} coincides with the $BMO(\mu)$, defined by \mref{bmodef}, if $\mu$ is doubling. It was only assumed in \cite{tolhh} that
\bdes \item[M.1)] There are a constant $C_\mu$ and a fixed number $n\in(0,d]$ such that for any cube $Q_r$ with side length $r>0$ \beqno{sizecond}\mu(Q_r)\le C_\mu r^n.\eeq \edes

The Hardy space $H^1(\mu)$ was introduced in \cite{tol3} and the duality $RBMO(\mu)$-$H^1(\mu)$ was also established. For our purpose in this paper, we don't need such a full force generality and we just recall the following deep result in \cite{tolhh}.

\blemm{rbmoreplem} (The Main Lemma - \cite[Lemma 4.1]{tolhh}) Let $f \in RBMO(\mu)$ with compact support and
$\iidmu{\Og}{f}=0$. There exist functions $h_m \in L^\infty(\mu)$ and $\fg_{y;m}$, $m \ge0$, such that
\beqno{rbmorep} f(x) = h_0(x) +\sum_{m=1}^\infty\int
\fg_{y;m}(x) h_m(y) d\mu(y);\eeq
with convergence in $L^1(\mu)$ and
\beqno{rbmorep1}\sum_{m=0}^{\infty}\|h_m\|_{L^\infty(\mu)}\le C[f]_{*,\mu}.\eeq 

Importantly, the functions $\fg_{y;m}$ satisfy the properties in \reflemm{fgproplem} below.
\elemm

It was shown in \cite{tolhh} that the functions $\fg_{y;m}$ satisfy the following properties.

\blemm{fgproplem} There is a constant $C$, depending also on $C_\mu$, such that for any $y\in \mbox{supp}(\mu)$ there is {\em some} cube $Q\subset \RR^d$ centered at $y$ \bdes \item[1)] $\fg_{y;m}\in C^1_0(Q)$. \item[2)] $0\le \fg_{y;m}(x)\preceq C l(Q)^{-n}$ for all $x\in Q$.
\item[3)] $|D\fg_{y;m}(x)|\preceq C l(Q)^{-n-1}$ or all $x\in Q$.
\edes
\elemm

\bproof In \cite[Lemma 7.8]{tolhh}, for suitable and fixed constants $\ag,\bg$ and  {\em some} cubes $Q_1,Q_2$ concentric with $Q$ and $\ag l(Q_1) \le l(Q)\le \bg l(Q_2)$, 1) comes from a) of \cite[Lemma 7.8]{tolhh} as $\fg_{y;m}=0$ outside $Q_2$. Similarly, 2) comes from  \cite[b) and c) of Lemma 7.8]{tolhh} if we note that $l(Q)\preceq|y-x|$ for $x\in Q_2\setminus Q_1$. Finally, 3) comes from \cite[d) of Lemma 7.8]{tolhh}.
\eproof

Right after the statement of \cite[Lemma 4.1]{tolhh}, there is a short proof of the fact that the $H^1(\mu)$ norm of a function is bounded by $\|f\|_{L^1(\mu)}+\|M_\Fg f\|_{L^1(\mu)}$ ($M_\Fg f$ is defined in \cite[Definition 1.1]{tolhh} which is generally larger then the one defined in \mref{MFgdef} below). For our purpose in this paper,  we need only estimate $\myprod{f,g}$ with $g\in RBMO(\mu)$. We then state the following lemma.

\blemm{fgestlem} Let $f\in RBMO(\mu)$ with the representation \mref{rbmorep}.
Let $F\in L^1(\mu)$ such that \beqno{fh1cond}\iidmu{\Og}{F}=0,\; M_{\hat{\Fg}}F\in L^1(\mu), \eeq where \beqno{MFgdef} M_{\hat{\Fg}}F(y) =\sup_{m\ge1}\int_{\Og}\fg_{y;m}(x)F(x)d\mu(x).\eeq

Then \beqno{fgest}|\myprod{F,f}|\le C(\|F\|_{L^1(\mu)}+\|M_{\hat{\Fg}} F\|_{L^1(\mu)})[f]_{*,\mu}.\eeq
\elemm

\bproof We repeat the argument right after the statement of \cite[Lemma 4.1]{tolhh}. From \mref{rbmorep}, we have
$$|\myprod{F,f}|\le \left|\int_\Og Fh_0d\mu\right| +\left|\sum_{m=1}^\infty\int_\Og\int_\Og
F\fg_{y;m}(x) h_m(y) d\mu(y)d\mu(x)\right|.$$

Since \beqno{keyfphi}\left|\int_\Og
F\fg_{y;m}(x) h_m(y) d\mu(y)\right|\le M_{\hat{\Fg}}F(x)\|h_m\|_{L^\infty(\mu)},\eeq by the definition \mref{MFgdef} of $M_{\hat{\Fg}} F$, we have $$|\myprod{F,f}|\le \|F\|_{L^1(\mu)}\|h_0\|_{L^\infty(\mu)} +\|M_{\hat{\Fg}} F\|_{L^1(\mu)}\sum_{m=1}^{\infty}\|h_m\|_{L^\infty(\mu)}.$$

By \mref{rbmorep1}, the above gives the lemma. \eproof

Inpired by \reflemm{fgproplem}, we introduce the following definition.
\bdefi{fgpropdef} A function $\fg\in C^1(\RR^d)$ is said to be in $\breve{\Fg}$ if  for  any $y\in\RR^d$ and {\em some} cube $Q\subset \RR^d$ centered at $y$ and the constant $C$ as in \reflemm{fgproplem} \bdes \item[f.1)]  $0\le \fg(x)\preceq Cl(Q)^{-n}$ for all $x\in Q$.
\item[f.2)] $\fg\in C^1_0(Q)$ and $|D\fg(x)|\preceq Cl(Q)^{-n-1}$ or all $x\in Q$.
\edes

For any $F\in L^1(\mu)$ we define
\beqno{Mfdefnorm1} M_{\breve{\Fg}}F(y) =\sup_{\fg\in\breve{\Fg}}\int_{\Og}\fg(x)f(x)d\mu(x)\in L^1(\mu),\eeq
\beqno{Mfdefnorm}\|F\|_{\breve{\Fg}}=\|F\|_{L^1(\mu)}+\|M_{\breve{\Fg}} F\|_{L^1(\mu)}.\eeq

\edefi

By \reflemm{fgproplem}, the functions $\fg_{y;m}$ belong to $\breve{\Fg}$ so that $M_{\hat{\Fg}}F(y)\le M_{\breve{\Fg}}F(y)$. We now have from \reflemm{fgestlem} the following result.  
\blemm{fgestlemz} Let $f\in BMO(\mu)$ and $F\in L^1(\mu)$ such that
\beqno{Mfdef} \iidmu{\Og}{F}=0.\eeq
Then \beqno{fgestz}|\myprod{F,f}|\le C\|F\|_{\breve{\Fg}}[f]_{*,\mu}.\eeq
\elemm

We will also use the definition of the  {\em centered}  Hardy-Littlewood maximal operator acting on functions $F\in L^1_{loc}(\mu)$
\beqno{maximal} M(F)(y) = \sup_\eg\{\iidmu{B_\eg(y)}{F(x)}\,:\, \eg>0 \mbox{ and } B_\eg(y)\subset\Og\}.\eeq

We also note here the Muckenhoupt theorem for non doubling measures. By \cite[Theorem 3.1]{OP}, we have that if $w$ is an $A_q(\mu)$ weight then for any $F\in L^q(\mu)$ with $q>1$
\beqno{HL1}\iidmu{\Og}{M(F)^q w} \le C(C_\mu,[w]_q)\iidmu{\Og}{F^q w}.\eeq
In particular, \beqno{HL}\iidmu{\Og}{M(F)^q} \le C_\mu\iidmu{\Og}{F^q}.\eeq

\section{The Main and Technical Inequality} \label{GNineqsec}\eqnoset
In this section we will establish a our main global weighted Gagliardo-Nirenberg  interpolation inequality. The main results stated in \refsec{mainres} are just consequences of this inequality.

Throughout this section we will always use the following notations and hypotheses. First, we repeat the condition M).

\bdes\item[M)] Let $d\mu=\og dx$ for some $\og\in C^1(\Og,\RR^+)$. Suppose that there are  a constant $C$ and a fixed number $n\in(0,d]$ such that : for any cube $Q_r$ with side length $r>0$ \beqno{sizecond1}\mu(Q_r)\le Cr^n.\eeq Furthermore, $\Og,\mu$ satisfies PS).\edes

The following assumptions slightly generalize A.1)-A.3) as we do not assume \mref{logcondszmain} in P.1). The assumptions P.2), W) are exactly A.2),A.3).

\bdes
\item[P.1)] Let $K:\mbox{dom}(K)\to\RR^m$ be a $C^1$ map on a domain $\mbox{dom}(K)\subset\RR^m$ such that $K_U(U)^{-1}$ exists for $U\in \mbox{dom}(K)$. Again, we use the notation \mref{DKdef}, $ \DKTmU=(K_U(U)^{-1})^T$, and assume that $\mccK_U\in L^\infty(\mbox{dom}K)$.

Let $\Fg,\LLg:\mbox{dom}(K)\to\RR^+$ be $C^1$ positive functions. Assume that \beqno{KUcond} |\mccK(U)|\preceq \LLg(U)\Fg^{-1}(U) \mbox{ for all $U\in \mbox{dom}(K)$}.\eeq

We also define the matrix 
\beqno{Pimatrix} \mccP(U):=\Fg^2(U)\LLg^{-1}(U)\DKTmU,\eeq
\item[P.2)] Let $U:\Og\to \mbox{dom}(K)$ be a $C^2$ vector-valued function. satisfying
\beqno{boundary}\myprod{\og\Fg^2(U)\DKTmU DU,\vec{\nu}}=0\eeq on $\partial\Og$, where $\vec{\nu}$ is the outward normal vector of $\partial\Og$.\edes

\bdes\item[W)] Let \beqno{myWdef}\myPi(x):=\LLg^{p+1}(U(x))\Fg^{-p}(U(x)) \mbox{ for $x\in \Og$}.\eeq Asume that $[\myPi^{\ag}]_{\bg+1}$ is finite for some $\ag>2/(p+2)$ and $\bg<p/(p+2)$. \edes

We recall the definitions \mref{Idefm}-\mref{Idef1zmm}

\beqno{Idef} I_1:=\iidmu{\Og}{\Fg^2(U)|DU|^{2p+2}},\;
I_2:=\iidmu{\Og}{\LLg^2(U)|DU|^{2p-2}|D^2U|^2},\eeq
\beqno{Idef1z} \myIbar:=\iidmu{\Og}{|\LLg_U(U)|^2|DU|^{2p+2}},\eeq \beqno{Idef2}\myIbreve:=\iidmu{\Og}{|D\og|^2\og^{-2}\LLg^2(U)|DU|^{2p}},\eeq and furthermore introduce \beqno{Idef1} \myIhat:=\iidmu{\Og}{(|\mccP_U(U)|\LLg(U)\Fg^{-1}(U))^2|DU|^{2p+2}}.\eeq

The main result of this section is the following theorem.

\btheo{GNlobal} Assume M), P.1)-P.2) and W).  
Suppose that  the integrals in \mref{Idef}-\mref{Idef2} are finite.
Then   there are constants $C,C([\myPi^{\ag}]_{\bg+1})$  for which \beqno{GNglobalest}I_1\le C\|K(U)\|_{BMO(\mu)}^2\left[I_2+\myIbar+C([\myPi^{\ag}]_{\bg+1}) [I_2+\myIhat+\myIbreve]\right].\eeq The constant $C$ depends on $C_{PS}, C_\mu$ and the constant $C$ in \refdef{fgpropdef}.
\etheo

The proof of this theorem will be divided into several lemmas. First of all,
let $W=K(U)$.  We then have $DU=K_U(U)^{-1}DW$ so that, from the definition of $\DKTmU=(K_U^{-1})^T$, $|DU|^2=\myprod{\DKTmU DU,DW}$. Hence, using  the definition of $\mccP(U)=\Fg^2(U)\LLg^{-1}(U)\DKTmU$ in \mref{Pimatrix}, we can write
\beqno{DWest}\barr{lll}I_1&=&\iidx{\Og}{\myprod{|DU|^{2p}\LLg(U)\Fg^2(U)\LLg^{-1}(U)\DKTmU DU,DW}\og}\\&=&\iidx{\Og}{\myprod{|DU|^{2p}\LLg(U)\og\mccP(U)DU,DW}}.\earr\eeq

Using the boundary assumption \mref{boundary}, $\myprod{\LLg(U)\og\mccP(U) DU,\vec{\nu}}=0$, and applying
integration by parts to the last integral, we have 
$$I_1=-\iidx{\Og}{\myprod{\Div(|DU|^{2p}\LLg(U)\og\mccP(U)DU),W}}.$$

Therefore, for $G:=\Div(|DU|^{2p}\LLg(U)\og\mccP(U)DU)\og^{-1}$
\beqno{ibp} I_1=  -\iidmu{\Og}{\myprod{G,W}}.\eeq 

From \mref{boundary} and integrations by parts again, we see that $$\iidmu{\Og}{G}=\iidx{\Og}{\Div(|DU|^{2p}\LLg(U)\og\mccP(U)DU)}=0.$$

We will establish bounds for $\|G\|_{L^1(\mu)},\|M_{\breve{\Fg}} G\|_{L^1(\mu)}$  and show that
\beqno{gh1esta} \|G\|_{\breve{\Fg}}  \le C\left[I_2^\frac12+\myIbar^\frac12 +C([\myPi^\ag]_{\bg+1}) [\myIhat^\frac12 + I_2^\frac12+\myIbreve^\frac12]\right]I_1^\frac12.\eeq 
Once this is proved, we obtain from \mref{ibp} and \mref{fgestz} of \reflemm{fgestlemz}, which is applicable here by M), that $I_1 \le C\|K(U)\|_{RBMO(\mu)}\|G\|_{\breve{\Fg}}$. As we are assuming that $\mu$ is doubling, $\|K(U)\|_{RBMO(\mu)}\sim \|K(U)\|_{BMO(\mu)}$. We then obtain
$$I_1 \le  C\|K(U)\|_{BMO(\mu)}
\left[I_2^\frac12+\myIbar^\frac12 + C([\myPi^\ag]_{\bg+1}) [I_2^\frac12+\myIhat^\frac12 + \myIbreve^\frac12]\right]I_1^\frac12,$$ which yields \mref{GNglobalest} via a simple use of Young's inequality. The proof is then complete.

To prove \mref{gh1esta}, we  first estimate $\|M_{\breve{\Fg}} G\|_{L^1(\mu)}$ and note that $$ M_{\breve{\Fg}} G=\sup_{\fg\in\breve{\Fg}}\left|\iidmu{\Og}{\fg G}\right|=\sup_{\fg\in\breve{\Fg}}\left|\iidx{\Og}{\fg g}\right|$$ where \beqno{gdef}g:=G\og=\Div(|DU|^{2p}\LLg(U)\og\mccP(U)DU).\eeq

Therefore, we need to establish that there are constants $C,C([\myPi^\ag]_{\bg+1})$  for which
\beqno{gh1est} \iidmu{\Og}{\sup_{\fg\in\breve{\Fg}} \left|\iidx{\Og}{\fg g}\right|} \le C\left[I_2^\frac12+\myIbar^\frac12+C([\myPi^\ag]_{\bg+1}) [\myIhat^\frac12 + I_2^\frac12+\myIbreve^\frac12]\right]I_1^\frac12.\eeq

From \mref{gdef} we can write $g=g_1+g_2$ with $g_i=\Div V_i$, setting 
\beqno{hJ0def} h:=\LLg(U) |DU|^{p-1}DU,\; J_{0,\eg}:=h_{B_\eg}=\mitmu{B_\eg}{\LLg(U) |DU|^{p-1}DU},\eeq
\beqno{g1def}V_1= \og|DU|^{p+1}\mccP(U)\left(h-J_{0,\eg}\right),\eeq 
\beqno{g2def}V_2= \og|DU|^{p+1}\mccP(U)J_{0,\eg}.\eeq 

We will establish \mref{gh1est} for $g$ being $g_1,g_2$ in the following lemmas. 

In the sequel, for any $\fg_\eg\in\breve{\Fg}$  and any $y\in\RR^d$ we denote by $B_\eg=B_\eg(y)$ the corresponding cube centered at $y$ with side length $\eg$ as in \refdef{fgpropdef}.

Let us consider $g_1$ first. 
\blemm{g1lemm} There is a constant $C$ such that  \beqno{g1est}\iidmu{\Og}{\sup_{\fg_\eg\in\breve{\Fg}} \left|\iidx{\Og}{\fg_\eg g_1}\right|} \le C
\left[I_2^\frac12+\myIbar^\frac12\right]I_1^\frac12.\eeq The constant $C$ depends on $C_{PS}, C_\mu$.
\elemm

\bproof   We use integration by parts (the boundary integral is zero because $\fg_\eg\in C^1_0(B_\eg)$)  to get
\beqno{gg1est}\barr{lll}\left|\iidx{B_\eg(y)}{\fg_\eg(x)g_1}\right|&=&\left|\iidmu{B_\eg(y)}{D\fg_\eg(x)\mccP(U)(h-J_{0,\eg})|DU|^{p+1}}\right|\\&\le&\frac{C}{\eg}
\mitmu{B_\eg(y)}{|h-h_{B_\eg(y)}||\mccP(U)||DU|^{p+1}}. \earr\eeq
Here,  we used the property of $D\fg_\eg$ in \refdef{fgpropdef}, which states $|D\fg_\eg|\preceq\eg^{-n-1}$, and the assumption M) that $\mu(B_\eg)\preceq\eg^{n}$. 

Note that \mref{KUcond} is equivalent to \beqno{piphi} |\mccP(U)|\preceq\Fg(U).\eeq
This and a simple use of H\"older's inequality for $q>2$ and \mref{piphi} yield that the last integral in \mref{gg1est} is bounded by
$$\frac{C}{\eg}
\left(\mitmu{B_\eg(y)}{|h-h_{B_\eg(y)}|^q}\right)^\frac1q
\left(\mitmu{B_\eg(y)}{[\Fg(U)|DU|^{p+1}]^{q'}}\right)^\frac1{q'}.
$$

Applying
the Poincar\'e-Sobolev inequality \mref{PSineq} to each component of $h$ and noting that there is a constant $C$ such that $$|Dh| \preceq |\LLg_U(U)||DU|^{p+1} + \LLg(U)|DU|^{p-1}|D^2U|,$$ we find a constant $C$ depends on $C_{PS}$ such that
\beqno{dhest}\barr{ll}\lefteqn{\frac{1}{\eg}
	\left(\mitmu{B_\eg}{|h-h_{B_\eg}|^q}\right)^\frac1q \le C
	\left(\mitmu{\tau_*B_\eg}{|Dh|^{q_*}}\right)^\frac1{q_*}}\hspace{.1cm}&\\& \le
C\left[\mitmu{\tau_*B_\eg}{(|\LLg_U(U)|^{q_*}|DU|^{(p+1)q_*}+
	\LLg^{q_*}(U)|DU|^{(p-1)q_*}|D^2U|^{q_*})}\right]^\frac1{q_*}.
\earr\eeq

Using the defintion of maximal functions \mref{maximal} and combining the above estimates, we get from \mref{gg1est}
\beqno{g1a}\sup_{\fg_\eg\in\breve{\Fg}} \left|\iidx{\Og}{\fg_\eg g_1}\right| \le C
\left[\Psi_1(y)+\Psi_2(y)\right]\Psi_3(y),\eeq
where $\Psi_i(y)=(M(F_i^{q_i}(y)))^\frac{1}{q_i}$ with $q_1=q_2=q_*$ and $q_3=q'$ and $$F_1=\LLg_U(U)|DU|^{p+1},\; F_2=\LLg(U)|DU|^{p-1}|D^2U|,\; F_3=\Fg(U)|DU|^{p+1}.$$

Because $q_i<2$ (as $q>2$ and $q_*=q\sg<2$), Muckenhoupt's inequality \mref{HL} implies
$$\left(\iidmu{\Og}{\Psi_i^2}\right)^\frac12=\left(\iidmu{\Og}{M(F_i^{q_i})^{\frac{2}{q_i}}}\right)^\frac12 \le  C_\mu\left(\iidmu{\Og}{F_i^{2}}\right)^\frac12.$$ 
Therefore, applying Holder's inequality to \mref{g1a} and using the above estimates and the notations \mref{Idef} and \mref{Idef1z}, we obtain \mref{g1est}. \eproof

\brem{PSmurem} We remark that \mref{dhest} is the only place where we need the assumption PS) that $\Og,\mu$ support a Poincar\'e-Sobolev inequality.  \erem

We now turn to $g_2$.

\blemm{g2estlemm} For any $p\ge1$ and $r\in(\frac{1}{p+1},1)$  we denote $$\ag(r)=\frac{r+1}{rp+r+1},\; \bg(r)=\frac{r(p+1)-1}{r(p+1)+1}.$$ 

Then there is $C([\myPi^{\ag(r)}]_{\bg(r)+1})\sim [\myPi^{\ag(r)}]_{\bg(r)+1}^\frac{1}{\ag(r)(p+1)}$ such that  \beqno{g2est}\iidmu{\Og}{\sup_{\fg_\eg\in\breve{\Fg}} \left|\iidx{\Og}{\fg_\eg g_2}\right|} \le C([\myPi^{\ag(r)}]_{\bg(r)+1})[\myIhat^\frac12 + I_2^\frac12+\myIbreve^\frac12 ]I_1^\frac12.\eeq 
\elemm

\bproof Note that  $\Div V_2 \le C(J_1+J_2+J_3)$ for some constant $C$ and $$J_1:= \og|\mccP_U||DU|^{p+2}J_{0,\eg}, \;  J_2:= \og|\mccP(U)||DU|^p|D^2U|J_{0,\eg},$$
$$J_3:=D\og|\mccP(U)||DU|^{p+1}J_{0,\eg},$$ with $J_{0,\eg}$ being defined in \mref{hJ0def}.

Because $\mccP(U):=\Fg^2(U)\LLg^{-1}(U)\DKTmU$ $$|\mccP_U|\preceq |(\Fg^2(U)\LLg^{-1})_U(U)||\mccK(U)|+\Fg^2(U)\LLg^{-1}(U)|\mccK_U(U)|.$$ We thus need only that $\mccK_U\in L^\infty(\mbox{dom}K)$ and so does $\mccP_U$. Our calculations for $J_1$ below are valid, see \cite[Theorem 7.8]{GT}.

In the sequel, for any $r>1/(p+1)$ we denote $r^*=1-\frac{1}{r(p+1)}$. We also write $f=\Fg|DU|^{p+1}$.

We consider $J_{0,\eg}$. From the notation $\myPi:=\LLg^{p+1}\Fg^{-p}$ (see \mref{myWdef})
$$J_{0,\eg}(y)\le\left|\mitmu{B_\eg}{\LLg\Fg^{\frac{-p}{(p+1)}}\Fg^{\frac{p}{p+1}}|DU|^p}\right|=\left|\mitmu{B_\eg}{\myPi^{\frac{1}{(p+1)}}f^{\frac{p}{p+1}}}\right|.$$

If $r_1>1/(p+1)$ we apply H\"older's inequality to the last integral to have the following estimate for $J_{0,\eg}$.
\beqno{J3est}J_{0,\eg}\le\left(\mitmu{B_\eg}{\myPi^{\frac{1}{r_1^*(p+1)}}}\right)^{r_1^*}\left(\mitmu{B_\eg}{f^{pr_1}}\right)^\frac{1}{r_1(p+1)}.\eeq

For $J_1$, we write $J_1=\og L_*LJ_{0,\eg}$ with $$L_*=|\mccP_U|\LLg\Fg^{-1}|DU|^{p+1},\; L=\LLg^{-1}\Fg|DU|.$$ By f.1) in \refdef{fgpropdef}, we have $\fg_\eg(x)\preceq \eg^{-n}\sim\mu(B_\eg)^{-1}$ so that we can use H\"older's inequality to get for any $s>1$
$$\sup_{\fg_\eg\in \breve{\Fg}} \left|\iidx{\Og}{\fg_\eg J_1}\right| \le \left(\mitmu{B_\eg}{L_*^{s'}}\right)^\frac{1}{s'}\left(\mitmu{B_\eg}{L^{s}}\right)^\frac{1}{s}J_{0,\eg}.$$

We write $L^s=\LLg^{-s}\Fg^{\frac{-sp}{(p+1)}}\Fg^{\frac{s}{p+1}}|DU|^s$
and use H\"older's inequality to have for any $r>1/(p+1)$ the following estimate.
$$\left(\mitmu{B_\eg}{L^s}\right)^\frac1s \le  \left(\mitmu{B_\eg}{|\LLg|^\frac{-s}{r_*}\Fg^{\frac{sp}{r^*(p+1)}}}\right)^{\frac{r^*}{s}}\left(\mitmu{B_\eg}{f^{sr}}\right)^\frac{1}{rs(p+1)}.$$

Combining these estimates with \mref{J3est} we then have \beqno{phiJ1est}\sup_{\fg_\eg\in \breve{\Fg}} \left|\iidx{\Og}{\fg_\eg J_1}\right| \le C_1M(L_*^{s'})^\frac{1}{s'}M(f^{sr})^\frac{1}{rs(p+1)}
M(f^{pr_1})^\frac{1}{r_1(p+1)},\eeq where, as $\myPi:=\LLg^{p+1}\Fg^{-p}$,
$$\barr{lll}C_1&\preceq&\left(\mitmu{B_\eg}{\myPi^{\frac{1}{r_1^*(p+1)}}}\right)^{r_1^*}\left(\mitmu{B_\eg}{|\LLg|^\frac{-s}{r_*}\Fg^{\frac{sp}{r^*(p+1)}}}\right)^{\frac{r^*}{s}}\\&=&\left[\left(\mitmu{B_\eg}{\myPi^{\frac{1}{r_1^*(p+1)}}}\right)
\left(\mitmu{B_\eg}{\myPi^{\frac{-s}{r^*(p+1)}}}\right)^{\frac{r^*}{r_1^*s}}\right]^{r_1^*}.\earr$$

We now choose $s,r,r_1$ such that $s'=sr=pr_1$ and $sr<2$. This is the case if $r\in(\frac{1}{p+1},1)$, $s=(r+1)/r$ then $s'=r+1$ and $r_1=(r+1)/p>1/(p+1)$. Let $\ag(r)=\frac{1}{r_1^*(p+1)}$ and $\bg(r)=\frac{r^*}{r_1^*s}$.   With such choice of $s,r,r_1$ we have \beqno{abchoice}\ag(r)=\frac{r+1}{rp+r+1},\; \bg(r)=\frac{r(p+1)-1}{r(p+1)+1},\; \ag(r)/\bg(r)=\frac{s}{r^*(p+1)}.\eeq 
It is clear that $C_1\preceq C_{1,r}^{r_1^*}$ with $r_1^*=\frac{1}{\ag(r)(p+1)}$ and
\beqno{C1def}
C_{1,r}=\sup_{B\subset\Og}\left(\mitmu{B}{\myPi^{\ag(r)}}\right)
\left(\mitmu{B}{\myPi^{\frac{-\ag(r)}{\bg(r)}}}\right)^{\bg(r)},\eeq where the supremum is taken over all cubes $B$ in $\Og$. Clearly, the definition of weight \mref{aweight} implies
\beqno{aweight-est}  [w]_{\nu+1}=\sup_{B\subset\Og}\left(\mitmu{B}{w}\right) \left(\mitmu{B}{w^{-\frac{1}{\nu}}}\right)^{\nu} \quad \mbox{for all $\nu>0$}.\eeq

From \mref{aweight-est}, $C_{1,r}=[\myPi^{\ag(r)}]_{\bg(r)+1}$. We then have \beqno{C1est} C_1\preceq C_{1,r}^{r_1^*}\preceq  [\myPi^{\ag(r)}]_{\bg(r)+1}^\frac{1}{\ag(r)(p+1)}.\eeq

As $sr=pr_1$, we then obtain from \mref{phiJ1est} the following.
$$\sup_{\fg_\eg\in\breve{\Fg}} \left|\iidx{\Og}{\fg_\eg J_1}\right| \preceq C_{1,r}^{r_1^*}
M(L_*^{sr})^\frac{1}{rs}M(f^{sr})^\frac{1}{rs}.$$ 

Applying H\"older's inequality to the right hand side, we get
$$\barr{ll}\lefteqn{\iidmu{\Og}{\sup_{\fg_\eg\in\breve{\Fg}} \left|\iidx{\Og}{\fg_\eg J_1}\right|} \preceq}\hspace{2cm}&\\& C_{1,r}^{r_1^*} \left(\iidmu{\Og}{M(L_*^{sr})^\frac{2}{rs}}\right)^\frac12\left(\iidmu{\Og}{M(f^{sr})^\frac{2}{rs}}\right)^\frac{1}{2}.\earr$$

Because $q=2/(rs)>1$, we can apply \mref{HL} to the integrals on the right and then use the definitions of $L_*,f, \myIhat$ to see that
\beqno{j1est}\iidmu{\Og}{\sup_{\fg_\eg\in\breve{\Fg}} \left|\iidx{\Og}{\fg_\eg J_1}\right|} \preceq C_{1,r}^{r_1^*}\|L_*\|_{2}\|f\|_{2}=C(C_{1,r})\myIhat^\frac12I_1^\frac12.\eeq

Next, we write $J_2=\og|\mccP(U)||DU|^{p-1}|D^2U||DU|J_{0,\eg}=\og L_*LJ_{0,\eg}$ with $$L_*=\LLg|DU|^{p-1}|D^2U|,\; L=\LLg^{-1}(U)|\mccP(U)||DU|.$$ 

We repeat the argument for $J_1$. Note that $|\mccP|\le\Fg$, by \mref{piphi}, and therefore $L^s\le\LLg^{-s}\Fg^{\frac{-sp}{(p+1)}}\Fg^{\frac{s}{p+1}}|DU|^s$. We have the following inequality.
$$\left(\mitmu{B_\eg}{L^s}\right)^\frac1s \le  \left(\mitmu{B_\eg}{|\LLg|^\frac{-s}{r_*}\Fg^{\frac{sp}{r^*(p+1)}}}\right)^{\frac{r^*}{s}}\left(\mitmu{B_\eg}{f^{sr}}\right)^\frac{1}{rs(p+1)}.$$ The estimate \mref{phiJ1est} for $J_1$ now applies to $J_2$ and yields \beqno{phiJ2est}\sup_{\fg_\eg\in\breve{\Fg}} \left|\iidx{\Og}{\fg_\eg J_2}\right| \le C_1M(L_*^{s'})^\frac{1}{s'}M(f^{sr})^\frac{1}{rs(p+1)}
M(f^{pr_1})^\frac{1}{r_1(p+1)}.\eeq 

As $sr=pr_1$, we  have as before $$\sup_{\fg_\eg} \left|\iidx{\Og}{\fg_\eg J_2}\right| \preceq C_{1,r}^{r_1^*}
M(L_*^{sr})^\frac{1}{rs}M(f^{sr})^\frac{1}{rs}.$$ 

The same argument for \mref{j1est}  for $J_1$  with the new definition of $L_*$ yields
\beqno{j2est}\iidmu{\Og}{\sup_{\fg_\eg\in\breve{\Fg}} \left|\iidx{\Og}{\fg_\eg J_2}\right|} \preceq C_{1,r}^{r_1^*}\|L_*\|_{2}\|f\|_{2}=C(C_{1,r})I_2^\frac12I_1^\frac12.\eeq

Concerning $J_3$, we write $J_3=D\og|\mccP(U)||DU|^{p}|DU|J_{0,\eg}=\og L_*LJ_{0,\eg}$ with $$L_*=D\og\og^{-1}\LLg|DU|^{p},\; L=\LLg^{-1}|\mccP(U)||DU|.$$ 
Similar argument for $J_2$ applying to this case then yields
\beqno{j4est}\iidmu{\Og}{\sup_{\fg_\eg\in\breve{\Fg}} \left|\iidx{\Og}{\fg_\eg J_3}\right|} \preceq C_{1,r}^{r_1^*}\|L_*\|_{2}\|f\|_{2}=C(C_{1,r})\myIbreve^\frac12I_1^\frac12\eeq

Combining the estimates  \mref{j1est},\mref{j2est} and \mref{j4est}, we derive \mref{g2est}. \eproof

Finally, we easily estimate $\|G\|_{L^1(\mu)}$. 

\blemm{divgL1} We have $$\iidmu{\Og}{|G|}\preceq CI_2^\frac12I_1^\frac12+C[\myIhat^\frac12 + I_2^\frac12+\myIbreve^\frac12 ]I_1^\frac12.$$\elemm
\bproof  Recall that $G:=g\og^{-1}$ so that $\|G\|_{L^1(\mu)}=\|g\|_{L^1(dx)}$. We write $g=\Div(B|DU|^{2p}DU)$ with $B=\og\LLg(U)\mccP(U)$. First of all, $$|\Div(B|DU|^{2p}DU)|\preceq |B_U||DU|^{2p+2}+|B||DU|^{2p}|D^2U|.$$ 

Because $|B_U|$ is bounded by a muliple of $$\{|D\og|\og^{-1}\LLg(U)|\mccP(U)|+|\LLg_U||\mccP(U)|+\LLg(U)|\mccP_U(U)|\}|DU|^{p+1}\LLg|Du|^p\og,$$ we see that a simple use of H\"older's inequality as in the proof of \reflemm{g2estlemm}, treating the last factor $\og\LLg|Du|^p$ as $J_0$, implies $$\iidx{\Og}{|B_U||DU|^{2p+2}}\le C[\myIhat^\frac12 + I_2^\frac12+\myIbreve^\frac12 ]I_1^\frac12.$$ 

As $|\mccP(U)|\preceq\Fg$, we have $|B||DU|^{2p}|D^2U|\preceq\Fg|DU|^{p+1}\LLg|DU|^{p-1}|D^2U|\og$. By H\"older's inequality we then obtain $$\iidx{\Og}{B|DU|^{2p+2}}\le CI_2^\frac12I_1^\frac12.$$ Combining the above estimates, we prove the lemma. \eproof

{\bf Proof of \reftheo{GNlobal}:} It is now clear that the above lemmas yield \beqno{ggestz}\|G\|_{\breve{\Fg}}\le C\left[I_2^\frac12+\myIbar^\frac12 + C([\myPi^{\ag(r)}]_{\bg(r)+1}) [\myIhat^\frac12 + I_2^\frac12+\myIbreve^\frac12]\right]I_1^\frac12.\eeq

Recall that
$\ag(r)=\frac{r+1}{rp+r+1}$ and $ \bg(r)=\frac{r(p+1)-1}{r(p+1)+1}$. We see that $\ag(r)$ decreases to $2/(p+2)$ and $\bg(r)$ increases to $p/(p+2)$ as $r\to1^-$.

From the definition of weights, a simple use of H\"older's inequality gives \beqno{wagbg} [w^\dg]_\cg\le [w]_\cg^\dg \quad \forall \dg\in(0,1).\eeq

Thus, if $\ag>2/(p+2)$ and $\bg<p/(p+2)$  then for $r$ close to $ 1$ we have $\ag(r)<\ag$ and $\bg(r)>\bg$.  
Hence, by choosing $r$ close to 1 and using \mref{wagbg} and the open end property of weights, we  see that
\beqno{C12est}[\myPi^{\ag(r)}]_{\bg(r)+1}\preceq C[\myPi^{\ag}]_{\bg+1}^\frac{\ag(r)}{\ag}.\eeq

Hence, we can replace $C([\myPi^{\ag(r)}]_{\bg(r)+1})$ by $C([\myPi^{\ag}]_{\bg+1})$ in \mref{ggestz}, which yields \mref{gh1est}. As $C([\myPi^{\ag(r)}]_{\bg(r)+1})\sim [\myPi^{\ag(r)}]_{\bg(r)+1}^\frac{1}{\ag(r)(p+1)}$, we can take $$C([\myPi^{\ag}]_{\bg+1})\sim [\myPi^{\ag}]_{\bg+1}^\frac{1}{\ag(p+1)}.$$

As we explain earlier, \mref{gh1est} yields $$I_1\le C\|K(U)\|_{BMO}\left[I_2^\frac12+\myIbar^\frac12 + C([\myPi^{\ag}]_{\bg+1}) [\myIhat^\frac12 + I_2^\frac12+\myIbreve^\frac12]\right]I_1^\frac12.$$
This gives
\beqno{GNglobalestkey}I_1\le C\|K(U)\|_{BMO}^2\left[I_2+\myIbar+C([\myPi^{\ag}]_{\bg+1}) [I_2+\myIhat+\myIbreve]\right].\eeq

The proof of the theorem is complete. 
\eproof

\brem{PSh} The only place we use the assumption PS) is \mref{dhest}. We just need to assume that PS) holds true for $h=\LLg(U)|DU|^{p-1}DU$ and some measure $\mu$ satisfying M). Combining with \cite[Theorem 5.1]{Haj} (see \refrem{PSrem}), which deals only with a pair $u,Du$, we need only that some Poincar\'e's inequality \mref{Pineq}, holds for the pair $h,Dh$. That is, we do not need \mref{Pineq} holds for any $h$ but the function $h=\LLg(U)|DU|^{p-1}DU$ in the consideration. \erem

\section{The Local Inequality} \label{localsec}\eqnoset

In this section, we will establish a local version of \reftheo{GNlobal}. Let $\Og_*$ be a subset of $\Og$. We assume that there are two functions $\og_*,\og_0$ satisfying the following conditions.

\bdes \item[L.0)] 
$\og_*\in C^1_0(\Og)$ and satisfies $\og_*\equiv 1$ in $\Og_*$ and $\og_*\le 1$ in $\Og$.
\beqno{subog} \og_*\equiv 1 \mbox{ in $\Og_*$ and } \og_*\le 1 \mbox{ in $\Og$}.\eeq

\item[L.1)] $\og_0\in C^1(\Og)$ and for $d\mu=\og_0^2dx$ and some $n\in(0,d]$ we have $\mu(B_r)\le Cr^n$. 

\item[L.2)]  The measure $\og_0^2 dx$ supports the Poincar\'e-Sobolev inequality \mref{PSineq} in PS). In addition, $\og_0$ also supports a Hardy type inequality: For any function $u\in C^1_0(B)$\beqno{lehr1} \iidx{\Og}{|u|^2|D\og_0|^2}\le C_H\iidx{\Og}{|Du|^2\og_0^2}.\eeq 
\edes

\btheo{GNlocalog1} Suppose P.1)-P.2) with $\og=\og_*\og_0^2$. Assume further that L.0)-L.2) hold true.
For any $\og_1\in L^1(\Og)$ and $\og_1\sim \og_0^2$ we define $d\mu=\og_1 dx$ and recall the definitions \mref{Idef}-\mref{Idef1} and introduce \beqno{I1*}I_{1,*}:=\iidmu{\Og_*}{\Fg^2(U)|DU|^{2p+2}},\eeq \beqno{I0*} \breve{I}_{0,*}:=\sup_\Og|D\og_*|^2\iidmu{\Og}{\LLg^2(U)|DU|^{2p}}.\eeq

Then, for any $\eg>0$ there are constants $C, C([\myPi^{\ag}]_{\bg+1})$ such that
\beqno{GNlocog11}I_{1,*}\le  \eg I_1+\eg^{-1}C\mC_*^2[I_2+\myIbar+C([\myPi^{\ag}]_{\bg+1}) [I_2+\myIhat+\myIbar+\breve{I}_{0,*}].\eeq
Here, $\mC_*:=\|K(U)\|_{BMO(\mu)}$ and $C$ depends on $C_{PS}$, $C_\mu$ and $C_H$.

\etheo

\bproof We consider first the case  $\og_1=\og_0^2$. Clearly, from the definition of $I_{1,*}$ and \mref{DWest}, we have for $W=K(U)$
$$I_{1,*}\le\int_{\Og}\Fg^2(U)|DU|^{2p+2}\og_*d\mu=\iidx{\Og}{\myprod{|DU|^{2p}\LLg(U)\og_*\og_0^2\mccP(U)DU,DW}}.$$

As we are assuming P.2) with $\og=\og_*\og_0^2$, \mref{boundary} gives \beqno{boundaryz}\myprod{\og_*\og_0^2\Fg^2(U) \DKTmU DU,\vec{\nu}}=0\eeq on $\partial\Og$ where $\vec{\nu}$ is the outward normal vector of $\partial\Og$.
Using this and integration by parts, we obtain 
$$I_{1,*}\le  -\iidmu{\Og}{\myprod{G,W}}$$ for $G:=\Div(|DU|^{2p}\LLg(U)\og_*\og_0^2\mccP(U)DU)\og_0^{-2}$ and $W=K(U)$. 

We now follow the proof of \reftheo{GNlobal} to establish a similar version of \mref{gh1esta}, with $d\mu=\og_0^2dx$, to complete the proof.  
First of all, L.1) implies M.1) so that \reflemm{fgestlemz} is applicable here. We see that \mref{gh1esta} holds true if \mref{gh1est} does.
We then need only establish a similar version of \mref{gh1est}.
Again, we can write $g=g_1+g_2$ with $g_i=\Div V_i$, setting 
$$V_1= \og_*\og_0^2|DU|^{p+1}\mccP(U)\left(h-J_{0,\eg}\right),
\;V_2= \og_*\og_0^2|DU|^{p+1}\mccP(U)J_{0,\eg},$$ where  $h:=\LLg(U) |DU|^{p-1}DU$ and \beqno{hJ0defz}  J_{0,\eg}:=h_{B_\eg}=\mitmu{B_\eg}{\LLg(U) |DU|^{p-1}DU}.\eeq We revisit the lemmas giving the proof of \mref{gh1est} and estimate \beqno{ogog1} \iidmu{\Og}{\sup_{\fg_\eg \in \breve{\Fg}} \left|\iidx{\Og}{\fg_\eg g_i}\right|}, \quad i=1,2.\eeq

Since $\og_*\le1$ we can discard it in the estimates for $g_1$ after the use of integration by parts \mref{gg1est} in the proof of \reflemm{g1lemm}. Because  the measure $\mu$  supports a Poincar\'e-Sobolev's inequality \mref{PSineq},  we can repeat the argument in the proof of \reflemm{g1lemm} to obtain the same estimate for the integral in \mref{ogog1} with $i=1$.  Similarly, we drop $\og_*$ in $J_i$'s, with the exception of $J_3$, in the proof of \reflemm{g2estlemm} to estimate the integral in \mref{ogog1} with $i=2$. Therefore,
$$\iidmu{\Og}{\sup_{\fg_\eg\in\breve{\Fg}} \left|\iidx{\Og}{\fg_\eg g}\right|} \le
\left[I_2^\frac12+\myIbar^\frac12+C([\myPi^{\ag}]_{\bg+1})[\myIhat^\frac12 + I_2^\frac12+\breve{I}_{*}^\frac12 ]\right]I_1^\frac12.$$
Here, the term $\breve{I}_{*}$, replacing $\myIbreve$ in \mref{g2est}, comes from the estimate for $J_3=D(\og_*\og_0^2)|\mccP(U)||DU|^{p}|DU|J_{0,\eg}$. In fact, we write $J_3=\og_0^2 L_*LJ_{0,\eg}$ for $L_*=D(\og_*\og_0^2)\og_0^{-2}\LLg|DU|^{p}$ and $L=\LLg^{-1}|\mccP(U)||DU|$. We obtain the following version of \mref{j4est} (with $C_{1,r}$ being replaced by $[\myPi^{\ag}]_{\bg+1}$)  $$\iidmu{\Og}{\sup_{\fg_\eg} \left|\iidx{\Og}{\fg_\eg J_3}\right|} \le C([\myPi^{\ag}]_{\bg+1})\breve{I}_{*}^\frac12I_1^\frac12,$$ with $\breve{I}_{*}=\|L_*\|_2^2$. That is, $$\breve{I}_{*}=\iidmu{\Og}{|D(\og_*\og_0^2)|^2\og_0^{-4}\LLg^2|DU|^{2p}}=\iidx{\Og}{|D(\og_*\og_0^2)|^2\og_0^{-2}\LLg^2|DU|^{2p}}.$$

Because $|D(\og_*\og_0^2)|^2\preceq |D\og_*|^2\og_0^4+\og_*^2|D\og_0|^2\og_0^2$, we have $$\breve{I}_{*}\preceq \iidx{\Og}{|D\og_*|^2\og_0^2\LLg^2|DU|^{2p}}+\iidx{\Og}{\og_*^2|D\og_0|^2\LLg^2|DU|^{2p}}.$$ The first integral on the right hand side is less than $\breve{I}_{0,*}$, defined by \mref{I0*}. Meanwhile, we apply the Hardy inequality \mref{lehr1} in L.2) to the second integral for $u=\og_*\LLg|DU|^{p}$, which belongs to $C^1_0(\Og)$, and note that (as $\og_*\le1$) $$|Du|^2\preceq \LLg^2|DU|^{2p-2}|D^2U|^2+|\LLg_U|^2|DU|^{2p}+|D\og_*|^2\LLg^2|DU|^{2p}.$$
We then have \beqno{hardyuse}\iidx{\Og}{\og_*^2|D\og_0|^2\LLg^2|DU|^{2p}}\le C\iidx{\Og}{|Du|^2\og_0^2}\preceq I_2+\myIbar+\breve{I}_{0,*}.\eeq

Thus, we get the following version of \mref{gh1est} $$\iidmu{\Og}{\sup_{\fg_\eg\in\breve{\Fg}} \left|\iidx{\Og}{\fg_\eg g}\right|} \le
C\left[I_2^\frac12+\myIbar^\frac12+C([\myPi^{\ag}]_{\bg+1})[\myIhat^\frac12 + I_2^\frac12+\breve{I}_{0,*}^\frac12 ]\right]I_1^\frac12.$$
The constant $C$ depends on $C_{PS},C_\mu$ and $C_H$.
Similarly, \reflemm{divgL1} gives a similar estimate for $\|G\|_{L^1(\mu)}$.  We then apply \reflemm{fgestlemz} as before and use Young's inequality to prove \mref{GNlocog11} for the case $d\mu=\og_0^2 dx$. 

Finally, if $\og_1\sim \og_0^2$ then the integrals  in \mref{GNlocog11} with respect to the two measures are comparable, because $D\og_1,D\og_0$ are not involved, so that \mref{GNlocog11} holds true as well. The proof is complete. \eproof

\brem{og*rem} For simplicity we assumed in L.0) that $\og_*\in C^1_0(\Og)$. More generally, we need only that $u=\og_*\LLg|DU|^{p}\in C^1_0(\Og)$ so that the Hardy inequality can apply in \mref{hardyuse}. \erem

\section{Proof of the Main Theorems and Further Generalizations}\label{collsec}\eqnoset

In this section, we present the proof of our main theorems.
To begin we will state the following theorem which is an immediate consequence of the main technical result \reftheo{GNlobal} and the definitions of the integrals in \mref{Idef}-\mref{Idef1}.

\btheo{GNlobalANS1a} Assume as in \reftheo{GNlobal}. Assume further that\beqno{LLgPicond} |\mccP_U|\LLg\Fg^{-1}\preceq \Fg.\eeq 
Then there are constants $C,C([\myPi^{\ag}]_{\bg+1})$  for which \beqno{GNglobalestz}I_1\le C\|K(U)\|_{BMO}^2\left[I_2+\myIbar+C([\myPi^{\ag}]_{\bg+1}) [I_2+I_1+\myIbreve]\right].\eeq 

In addition, if \beqno{KeyLLgPicond1}|\LLg_U|\preceq \Fg\eeq
then \beqno{GNglobalesta}I_1\le C\|K(U)\|_{BMO}^2\left[I_2+I_1+C([\myPi^{\ag}]_{\bg+1}) [I_2+I_1+\myIbreve]\right].\eeq

\etheo

\bproof By \mref{LLgPicond} and the definition of $\myIhat$ in \mref{Idef1}, we have  $\myIhat\preceq I_1$.  Similarly, \mref{KeyLLgPicond1} and \mref{Idef1z}  give $\myIbar\le I_1$. This theorem then follows from \reftheo{GNlobal}.
\eproof

The local version \reftheo{GNlocalog1} then implies the following \bcoro{GNlocalzzz} Assume \mref{LLgPicond} and \mref{KeyLLgPicond1}. Using the definitions \mref{I1*} and \mref{I0*} for $I_{1,*}$ and $\breve{I}_{0,*}$, we have \beqno{GNlocog11ANSz}I_{1,*}\le  \eg I_1+\eg^{-1}C\|K(U)\|_{BMO(\Og)}^2[I_2+I_1+C([\myPi^{\ag}]_{\bg+1}) [I_2+I_1+\breve{I}_{0,*}]].\eeq
\ecoro

Concerning the condition \mref{LLgPicond}a and for later references, we remark the following.

\brem{LLgPicondrem} The technical theorems \reftheo{GNlobal} and \reftheo{GNlocalog1} always assume that (recalling $\DKTmU=(K_U^{-1}(U))^T$ and $|\DKTmU|=|K_U^{-1}(U)|$)
\beqno{mycond1}|\DKTmU|\preceq \LLg(U)\Fg^{-1}(U).\eeq 

The condition \mref{LLgPicond} can be replaced by a stronger but more verifiable one: \beqno{mycond2} |\Fg_U(U)||\DKTmU|\preceq\Fg(U),\; |\mccK_U(U)| \mbox{ is bounded}.\eeq

Indeed, from the definition $\mccP(U)=\Fg^2(U)\LLg^{-1}(U)\DKTmU$, we have $$|\mccP_U|\LLg\Fg^{-1} \preceq |\Fg_U||\mccK|+\Fg|\LLg_U|\LLg^{-1}|\mccK|+\Fg|\mccK_U(U)|.$$ By \mref{mycond1}, we have $\Fg|\LLg_U|\LLg^{-1}|\mccK|\preceq|\LLg_U|$. Thus,  if \mref{mycond2} holds then the above clearly implies $|\mccP_U|\LLg\Fg^{-1}\preceq \Fg+|\LLg_U|$ so that $\myIhat\preceq I_1+\myIbar$. Therefore, \mref{GNglobalestz} also holds true if \mref{mycond1} and \mref{mycond2} are assumed.\erem

{\bf Proof of \reftheo{GNlobalANSm} and \reftheo{GNlocalog1m}:}
The assumption A.2) contains \mref{mycond1} and \mref{mycond2} of \refrem{LLgPicondrem}. Therefore,  under the assumptions A.1)-A.3), \reftheo{GNlobalANSm} follows from \reftheo{GNlobalANS1a}. In the same way, \reftheo{GNlocalog1m} is a consequence of \refcoro{GNlocalog1}.
\eproof

Let us consider different choices of $K$ and prove \refcoro{GNlobalANS}.
Consider  the case $\LLg(U)=(\eg+|U|)^k$ and $\Fg(U)\sim|\LLg_U(U)|$ for any $k\ne0$ and $\eg\ge0$. The corresponding integrals in \reftheo{GNlobalANS1a} are
\beqno{logI1}I_1= \iidmu{\Og}{(\eg+|U|)^{2k-2}|DU|^{2p+2}},\;\myIbreve= \iidmu{\Og}{(\eg+|U|)^{2k}|DU|^{2p}},\eeq
\beqno{logI2}I_2= \iidmu{\Og}{(\eg+|U|)^{2k}|DU|^{2p-2}|D^2U|^2}.\eeq

{\bf Proof of \refcoro{GNlobalANS}:} We apply \mref{GNglobalesta} of \reftheo{GNlobalANS1a} to this case.  $\LLg(U)=(\eg+|U|)^k$, $\Fg(U)=|k|(\eg+|U|)^{k-1}\sim |\LLg_U(U)|$. It is clear that $\myPi = \LLg^{p+1}\Fg^{-p}\sim (\eg+|U|)^{k+p}$.

As in \mref{Klogm}, we define $K(U)=[\log(\eg+|U_i|)]_{i=1}^m$, therefore $K_U(U)=\mbox{diag}[(\eg+|U_i|)^{-1}]$ and $\DKTmU=\mbox{diag}[(\eg+|U_i|)]$. Hence, $|\DKTmU|\preceq\eg+|U|$ so that $|\mccK|\preceq\LLg\Fg^{-1}$. Also, it is clear that $|\Fg_U||\mccK|\preceq\Fg$ and $|\mccK_U(U)|$ is bounded. Hence, \mref{GNglobalesta} holds true by \refrem{LLgPicondrem} and applies here to give $$I_1\le C\|K(U)\|_{BMO(\mu)}^2\left[I_2+I_1+C([\myPi^{\ag}]_{\bg+1}) [I_2+I_1+\myIbreve]\right],$$ and the proof is complete. \eproof

To prove \refcoro{GNlobalANS1} we have the following estimate for $[\myPi^\ag]_{\bg+1}$.

\blemm{Westlemm} For any $\ag,\bg>0$ and $\mmc_1,\mmc_2$ as in \mref{JNineq}
\beqno{Westcond} [\log(\myPi)]_{*,\mu}\le \mmc_1\bg\ag^{-1} \Rightarrow [\myPi^\ag]_{\bg+1}\le \mmc_2^{1+\bg}.\eeq
\elemm

\bproof We again recall the John-Nirenberg inequality \mref{JNineqm}: If $\mu$ is doubling then for any BMO($\mu$) function $v$ 
there are constants $\mmc_1,\mmc_2$, which depend only on the doubling constant of $\mu$, such that  
\beqno{JNineq}\mitmu{B}{e^{\frac{\mmc_1}{[v]_{*,\mu}}|v-v_B|}} \le \mmc_2.\eeq

For any $\bg>0$ we know that $e^v$ is an $A_{\bg+1}$ weight with $[e^v]_{\bg+1}\le \mmc_2^{1+\bg}$ (e.g. see \cite[Chapter 9]{Graf}) if  \beqno{Acgcond}\sup_B \mitmu{B}{e^{(v-v_B)}} \le \mmc_2,\; \sup_B \mitmu{B}{e^{-\frac{1}{\bg}(v-v_B)}} \le \mmc_2.\eeq

It is clear that \mref{Acgcond} follows from \mref{JNineq} if $\mmc_1[v]_{*,\mu}^{-1}\ge \max\{1,\bg^{-1}\}$.
Using these facts  with $v=\ag \log\myPi$, we see that \mref{Westcond} holds. \eproof

{\bf Proof of \refcoro{GNlobalANS1}:} From the definition of $\myPi=\LLg^{p+1}\Fg^{-p}=|k|^{-p}(\eg+|U|)^{k+p}$.
We then have $[\log(\myPi)]_{*,\mu}=|k+p|[\log(\eg+|U|)]_{*,\mu}$. Therefore the assumption \mref{Westcond1}, that $|k+p|[\log(\eg +|U|)]_{*,\mu}\le \mmc_1\bg\ag^{-1}$, and \mref{Westcond} imply  $[\myPi^\ag]_{\bg+1}\le \mmc_2^{1+\bg}$. The Corollary then follows from \refcoro{GNlobalANS}. \eproof

\end{document}